\def\L{{\mathcal L}}
\def\C{{\mathcal C}}
\def\A{{\mathcal A}}
\def\B{{\mathcal B}}
\def\N{\mathbb N}

\def\Q{\mathbb Q}
\def\pfz{\begin{proof}}
\def\pfk{\end{proof}}
\documentclass[10pt,intlimits]{amsart}
\usepackage[T1]{fontenc}
\usepackage[latin1]{inputenc}

\usepackage{amsthm,amsmath}
\usepackage{amssymb,color}

\usepackage{enumerate}
\usepackage{epic}

\newtheorem{lem}{Lemma}
\newtheorem{thm}[lem]{Theorem}
\newtheorem{prop}[lem]{Proposition}
\newtheorem{coro}[lem]{Corollary}
\newtheorem{pozn}[lem]{Remark}

\newcommand{\Hm}[1]{\leavevmode{\marginpar{\tiny%
$\hbox to 0mm{\hspace*{-0.5mm}$\leftarrow$\hss}%
\vcenter{\vrule depth 0.1mm height 0.1mm width \the\marginparwidth}%
\hbox to
0mm{\hss$\rightarrow$\hspace*{-0.5mm}}$\\\relax\raggedright #1}}}

\begin{document}
%%%%%%%%%%%%%%%%%%%%%%%%%%%%%%%%%%%%%%%%%%%%%%%%%%%%%%%%%%%%%%%%%%%%%
\title[Powers in 3iets]{Note on powers in three interval exchange transformations}
%-------------------------------------------------------------%
\author[D.~Lenz]{Daniel Lenz$^1$}
\address{$^1$ Mathematisches Institut, Friedrich-Schiller Universit\"at Jena,
  Ernst-Abb\'{e} Platz 2, 07743 Jena, Germany}
\email{ daniel.lenz@uni-jena.de }
\author[Z.~Mas\'{a}kov\'{a}]{Zuzana Mas\'{a}kov\'{a}$^2$}
\address{$^2$ Department of Mathematics and Doppler Institute, FNSPE, Czech Technical University in Prague, Trojanova 13, 120 00 Praha 2, Czech Republic }
\email{ zuzana.masakova@fjfi.cvut.cz }
\author[E.~Pelantov\'{a}]{Edita Pelantov\'{a}$^3$}
\address{$^3$ Department of Mathematics and Doppler Institute, FNSPE, Czech Technical University in Prague, Trojanova 13, 120 00 Praha 2, Czech Republic }
\email{ edita.pelantova@fjfi.cvut.cz }

\begin{abstract}
We study repetitions in infinite words coding exchange of three intervals with permutation (3,2,1), called 3iet words. The language of such words is determined by two parameters $\varepsilon,\ell$. We show  that  finiteness of the index of 3iet words is equivalent to  boundedness  of the  coefficients of the continued fraction of  $\varepsilon$. In this case we also  give an upper and lower estimate on the index of the corresponding 3iet word. \end{abstract}
%-------------------------------------------------------------%
\maketitle

%%%%%%%%%%%%%%%%%%%%%%%%%%%%%%%%%%%%%%%%%%%%%%%%%%%%%%%%%%%%%%%%%%%%%
\section{Introduction}

In this paper we study the index of infinite words of a certain class,
namely infinite words coding exchange of three intervals, the so-called 3iet words. The index of an infinite word $u$
describes maximal repetition of a factor of $u$. The study of repetitions in infinite words has become more extensive in recent years, also because of possible applications in the spectral theory of associated Schroedinger operators (see Remark (v) in Section \ref{sec:relationtoSturmanndRotation}).  Repetitions in the most prominent infinite words, namely Sturmian words, were studied by many authors. The first characterization of Sturmian words with finite index was given by Mignosi in~\cite{Mignosi}.

\begin{thm}\label{thm:1}
Let $\varepsilon\in(0,1)$ be an irrational number with continued fraction expansion of the form $[0,a_1,a_2,\dots]$. A Sturmian word $u_\varepsilon$
with the slope $\varepsilon$ has a finite index if and only if the sequence $(a_n)_{n=1}^\infty$ is bounded.
\end{thm}

The characterization of words with finite index is known also for another class of infinite words, namely for the words coding rotations. The language of these words depends on a pair of parameters, $(\alpha,\beta)$. Adamczewski in~\cite{adamczewski-jntb} introduces the notion of ${\mathcal D}$-expansion of such a pair. In~\cite{adamczewski-powers} he  then uses this notion to characterize words coding rotations which have finite index.

Our aim is to characterize 3iet words with finite index. This is the content of our main result,  Theorem~\ref{thm:XY},  which we prove in Section~\ref{sec:3ietpowers}. Before that, we introduce necessary notions and results, (Section~\ref{sec:preliminaries}). Although 3iet words are ternary words, they are closely connected to both mentioned classes of binary words: Every Sturmian word and every coding of rotation is an image of a 3iet word under a suitable morphism. In Section~\ref{sec:relationtoSturmanndRotation} we discuss the relation of results obtained for 3iet words with results known for these binary words.

%%%%%%%%%%%%%%%%%%%%%%%%%%%%%%%%%%%%%%%%%%%%%%%%%%%%%%%%%%%%%%%%%%%%%
\section{Preliminaries}\label{sec:preliminaries}

%--------------------------------------------------------------------
\subsection{Basic notions of combinatorics on words}

An alphabet $\mathcal A$ is a finite set of symbols, called
letters. A word $w$ of length $|w|=n$ is a concatenation of $n$ letters.
${\mathcal A}^*$ is the set of all finite words over the
alphabet $\A$ including the empty word $\epsilon$. Equipped with
the operation of concatenation, it is a monoid. We define also
infinite words $u=(u_n)_{n\in\N}\in\A^{\N}$.
A finite word $v\in{\A}^*$ is called a {\em factor} of a word $w$
(finite or infinite), if there exist words $w^{(1)},w^{(2)}$ such
that $w=w^{(1)}vw^{(2)}$. If $w^{(1)}=\epsilon$, then $v$ is said
to be a prefix of $w$, if $w^{(2)}=\epsilon$, then $v$ is a suffix
of $w$. The set of all factors of $u$ is
called the {\em language} of $u$ and denoted by $\L(u)$.
The occurrence of a  factor $w$ of an infinite word $u=u_0u_1u_2u_3\cdots$
 is an index $i\in\N$ such that $w$ is a prefix of the infinite word $u_iu_{i+1}u_{i+2}\cdots$.
An infinite word $u$ is said to be recurrent if every factor $w\in\L(u)$ has in $u$
at least two occurrences (and thus infinitely many occurrences).
If moreover the distances between consecutive occurrences of a given factor $w$ are bounded for every $w\in\L(u)$, then
the infinite word $u$ is called uniformly recurrent.

In this paper we focus on repetition of factors in infinite words. We
say that a word $v$ is a {\em power} of a word $w$, if
$|v|\geq|w|$ and $v$ is a prefix of the periodic word $www\cdots$.
We write $v=w^r$ where $r=|v|/|w|$. The index of a word $w$ in an
infinite word $u$ is defined by
\begin{equation}\label{e:indfactor}
{\rm ind}(w)= \sup\{r\in\Q\mid w^r\in\L(u)\}\,.
\end{equation}
It is easy to show that in uniformly recurrent non-periodic words, every factor $w$ has finite index.
Taking supremum of indices over all factors of an infinite word
$u$, one obtains an important characteristics of $u$, the
so-called {\em index} of $u$. Formally,
\begin{equation}\label{e:indword}
{\rm ind}(u)= \sup\{{\rm ind}(w)\mid w\in\L(u)\}\,.
\end{equation}

%--------------------------------------------------------------------
\subsection{Words coding exchange of three intervals}

We consider the index of infinite words coding exchange of three intervals with permutation $(3,2,1)$.
This transformation $T_{\varepsilon,\ell}$ depends on two parameters $\varepsilon,\ell$
which satisfy the condition
\begin{equation}\label{eq:3}
\varepsilon\in(0,1)\setminus\Q\,,\qquad
\max\{\varepsilon,1-\varepsilon\}<\ell<1\,.
\end{equation}
It is defined as  $T_{\varepsilon,\ell}(x):[0,\ell)\mapsto [0,\ell)$,
$$%\begin{equation}\label{eq:1}
T_{\varepsilon,\ell}(x):=
\begin{cases}
x+1-\varepsilon & \text{for } x\in I_A:=[0,\ell-1+\varepsilon)\,,\\
x+1-2\varepsilon &\text{for } x\in I_B:=[\ell-1+\varepsilon,\varepsilon)\,,\\
x-\varepsilon &\text{for } x\in I_C:=[\varepsilon,\ell)\,.
\end{cases}
$$%\end{equation}
The orbit
of a point $x_0\in[0,\ell)$ under the transformation
$T_{\varepsilon,\ell}$ can be coded by the infinite
word $u_{\varepsilon,\ell,x_0}=(u_n)_{n\in\N}$ in the alphabet
$\{A,B,C\}$, where
\begin{equation}\label{eq:4}
u_n=X\,,\qquad\hbox{if }\ T_{\varepsilon,\ell}^n(x_0)\in I_X\,.
\end{equation}
The infinite word $u_{\varepsilon,\ell,x_0}$ is called a 3iet word. Many combinatorial properties
of such words can be found in~\cite{ferenczi-holton-zamboni-jam-89}. An  equivalent approach to such words using
cut-and-project sequences is given in~\cite{gmp-jtnb-15}. Let us point out that a 3iet word is
aperiodic, uniformly recurrent, and its language (and inparticular also its index) does not depend on the
choice of the initial point $x_0$.
Consequently we will often omit $x_0$ in the notation of 3iet words.
%Consequently, for the language, resp. index of $u_{\varepsilon,\ell,x_0}$
%we use the notation $\L(u_{\varepsilon,\ell})$, resp. ${\rm ind}(u_{\varepsilon,\ell})$.

%--------------------------------------------------------------------
\subsection{Sturmian words}

3iet words can be considered as generalizations of Sturmian words. These can be defined by many equivalent ways, one of them
uses exchange of two intervals. For an irrational $\varepsilon\in(0,1)$ we define
the mapping $T_\varepsilon:[0,1)\mapsto[0,1)$ by
\begin{equation}\label{eq:1}
T_{\varepsilon}(x)=
\begin{cases}
x+1-\varepsilon &\text{for } x\in I_0:=[0,\varepsilon)\,,\\
x-\varepsilon   &\text{for } x\in I_1:=[\varepsilon,1)\,.
\end{cases}
\end{equation}
The infinite word $u_{\varepsilon,x_0}=(u_n)_{n\in\N}\in\{0,1\}^{\N}$ coding the orbit of a point $x_0\in [0,1)$ under $T_\varepsilon$
is given by
\begin{equation}\label{eq:2ietslovo}
u_n=
\begin{cases}
0 & \text{if } T_{\varepsilon}^n(x_0)\in I_0\,,\\
1 & \text{if } T_{\varepsilon}^n(x_0)\in I_1\,.
\end{cases}
\end{equation}
The word $u_{\varepsilon,x_0}$ is the upper mechanical word corresponding to
$\alpha=1-\varepsilon$ and $\beta=x_0$. Similarly, lower
mechanical words are obtained by coding exchange of intervals of
the type $(\cdot,\cdot]$. It is the representation of Sturmian
words by two interval exchange which will be the most used here.

The language and thus also the index of a Sturmian word $u_{\varepsilon,x_0}$ does not depend on the
initial point $x_0$, but rather only on the parameter
$\varepsilon$, which is called the slope of $u_{\varepsilon,x_0}$.
Therefore we use the notation $u_\varepsilon $ without specification of $x_0$.
%Therefore we use the notation
%${\mathcal L}(u_\varepsilon)$ and ${\rm ind}(u_\varepsilon)$ without specifiaction of $x_0$.

In case that the index of a Sturmian word is finite, its value is also described. The exact result has been independently obtained in~\cite{Carpi}
and~\cite{DaLenzNecelo}.

\begin{thm}\label{thm:3}
Let $\varepsilon\in(0,1)$ be an irrational number with continued fraction expansion of the form $[0,a_1,a_2,\dots]$. Then for a Sturmian word $u_\varepsilon$
with the slope $\varepsilon$ we have
$$
{\rm ind}(u_\varepsilon)=\sup\Big\{\,2+a_{N+1}+\frac{q_{N-1}-2}{q_N}\;\Big|\; N\geq
0\,\Big\}\,,
$$
where $q_N$ are the denominators of the convergents of $\varepsilon$.
\end{thm}

This theorem will be useful for the description of the index of 3iet words,  since
one can characterize 3iet words by Sturmian words, as has been shown in \cite{ABMP}. For the precise formulation of the result, we need the notion of a morphism $\varphi$
of the monoid $\A^*$ to the monoid $\B^*$, where $\A$, $\B$ are finite, in general different, alphabets.
Let us recall that a morphism $\varphi:\A^*\mapsto\B^*$ is a mapping satisfying $\varphi(vw)=\varphi(v)\varphi(w)$ for every $v,w\in\A^*$.
It is obvious that a morphism is uniquely determined by the values $\varphi(a)$ for every letter $a$ in the alphabet $\A$. The action of
a morphism can be naturally extended to an infinite word $u\in\A^{\N}$ by
$$
\varphi(u_0u_1u_2\cdots):=\varphi(u_0)\varphi(u_1)\varphi(u_2)\cdots\,.
$$

\begin{thm}\label{thm:ABMP}
Let $\sigma_{01},\sigma_{10}:\{A,B,C\}^*\mapsto\{0,1\}^*$ be morphisms defined by
\begin{equation}\label{eq:morfismysigma}
\begin{split}
\sigma_{01}:&\quad A \mapsto 0\,,\ B  \mapsto 01\,,\  C \mapsto 1\,,\\
\sigma_{10}:&\quad A \mapsto 0\,,\ B  \mapsto 10\,,\  C \mapsto 1\,.
\end{split}
\end{equation}
\begin{itemize}
\item
A ternary word $u\in\{A,B,C\}^{\N}$ is a 3iet word if and only if both
$\sigma_{01}(u)$ and $\sigma_{10}(u)$ are Sturmian words.
\item
If $u$ is a 3iet word with parameters $\varepsilon,\ell$, then $\sigma_{01}(u)$
and $\sigma_{10}(u)$ are Sturmian words with slope $\varepsilon$.
\end{itemize}
\end{thm}

The above theorem motivates the  notion of amicability of words discussed next.
We say that the word $w^{(1)}\in \{0,1\}^*\cup \{0,1\}^\mathbb{N}$ is {\em amicable} with $w^{(2)}\in \{0,1\}^*\cup \{0,1\}^\mathbb{N}$, if
there exists a word $w\in\{A,B,C\}^*\cup \{A,B,C\}^\mathbb{N}$ such that $w^{(1)}=\sigma_{01}(w)$ and
$w^{(2)}=\sigma_{10}(w)$, we then write  $w^{(1)}\propto w^{(2)}$. Note that this relation is not symmetric.
The word $w$ is called the {\em ternarization} of the pair $w^{(1)},w^{(2)}$
and denoted by $w={\rm ter}(w^{(1)},w^{(2)})$. This notion is illustrated on Figure~\ref{fff}.

\begin{figure}[!ht]
\begin{equation*}
\begin{aligned}
w^{(1)} =\quad \\[-0.9mm] w^{(2)} =\quad \\[1mm] w =\quad
\end{aligned}
\begin{gathered}\boxed{\:\begin{gathered}0 \\ 0\end{gathered}\:}\\A\end{gathered}\
\begin{gathered}\boxed{\:\begin{gathered}1 \\ 1\end{gathered}\:}\\C\end{gathered}\
\begin{gathered}\boxed{\:\begin{gathered}0 \\ 0\end{gathered}\:}\\A\end{gathered}\
\begin{gathered}\boxed{\:\begin{gathered}0\ 1 \\ 1\ 0\end{gathered}\:}\\B\end{gathered}\
\begin{gathered}\boxed{\:\begin{gathered}0 \\ 0\end{gathered}\:}\\A\end{gathered}\
\begin{gathered}\boxed{\:\begin{gathered}1 \\ 1\end{gathered}\:}\\C\end{gathered}
\end{equation*}
\caption{The word $w^{(1)}=0100101$ is amicable to $w^{(2)}=0101001$. Their ternarization is the word
$w=ACABAC$.} \label{fff}
\end{figure}

\begin{pozn} \label{uniqueness}
Note that the ternarization of $w^{(1)}$ and $w^{(2)}$ with $w^{(1)}$ amicable to $w^{(2)}$ is unique, as is clear from Figure~\ref{fff}.
\end{pozn}

%--------------------------------------------------------------------
\subsection{Words coding rotations}

The transformation $T_{\varepsilon}(x)$ of~\eqref{eq:1} defining Sturmian words can be rewritten in a more compact form
using the fractional part $\{y\}:=y-\lfloor y\rfloor$ as
$$
T_{\varepsilon}(x)=\{x+(1-\varepsilon)\}\,,\qquad\hbox{for }\ x\in[0,1)\,.
$$
In these terms, the $n$-th iteration of $T_{\varepsilon}(x)$ can be given as $T_{\varepsilon}^n(x)=\{x+n(1-\varepsilon)\}$,
and therefore the definition of the Sturmian word coding the orbit of $0$ can be rewritten simply as
$$%\begin{equation}\label{eq:2ietslovo}
u_n=
\begin{cases}
0 & \text{if } \{n(1-\varepsilon)\}\in [0,\varepsilon)=I_0\,,\\
1 & \text{if } \{n(1-\varepsilon)\}\in [\varepsilon,1)=I_1\,.
\end{cases}
$$%\end{equation}
Relaxing the relation between the length of the intervals $I_0,I_1$ and the slope $\varepsilon$, we obtain a generalization
of Sturmian words called coding of rotations. Such a word $u=(u_n)_{n\in\N}$ depends on two prameters $\alpha,\beta\in[0,1)$, $\alpha\notin\Q$,
and is given by
$$
u_n=
\begin{cases}
0 & \text{if } \{n\alpha\}\in [0,\beta)=I_0\,,\\
1 & \text{if } \{n\alpha\}\in [\beta,1)=I_1\,.
\end{cases}
$$
It is interesting to mention the relation of codings of rotations and 3iet words which was described in~\cite{adamczewski-jntb}: For every word $u\in\{0,1\}^{\N}$ coding of rotation, there exists a 3iet word $v\in\{A,B,C\}^{\N}$ and an integer $k\in\N$ such that $u=\Phi_k(v)$, where
$\Phi_k: \{A,B,C\}^*\mapsto \{0,1\}^*$ is a morphism
\begin{eqnarray*}
\Phi_k(A)&=&0\,,\\
\Phi_k(B)&=&01^{k+1}\,,\\
\Phi_k(C)&=&01^k\,.
\end{eqnarray*}
%Recall that for a finite word $w$, the symbol $w^k$ denoted $k$-fold repetition of the word $w$, $w^k=\underbrace{ww\cdots w}_{\hbox{$k$-times}}$.
One can give also the opposite statement: The infinite word $\Phi_k(v)$ is a coding of rotation for every 3iet word $v$ and every $k\in\N$.

%%%%%%%%%%%%%%%%%%%%%%%%%%%%%%%%%%%%%%%%%%%%%%%%%%%%%%%%%%%%%%%%%%%%%
\section{Powers in 3iet words}\label{sec:3ietpowers}

The aim of this section is a characterization of 3iet words with finite index. We prove the following theorem.

\begin{thm}\label{thm:XY}
Let $\varepsilon,\ell$ be parameters satisfying~\eqref{eq:3} and let $\varepsilon=[0,a_1,a_2,\dots]$ be the continued fraction expansion of $\varepsilon$.
Let $u\in\{A,B,C\}^{\N}$ be a 3iet word with parameters ${\varepsilon,\ell}$. Then
$$
{\rm ind}(u)<+\infty \qquad\iff\qquad \sup_{n\in\N} a_n  < +\infty\,.
$$
\end{thm}

The characterization of 3iet words having finite index is thus analogous to that proved in 1989 by Mignosi for Sturmian words, see~\cite{Mignosi}.
For the proof, we divide Theorem~\ref{thm:XY} into two implications.

\begin{prop}\label{p:1}
Let $u_{\varepsilon,\ell}$ be a 3iet word with parameters ${\varepsilon,\ell}$ and let ${\rm ind}(u_{\varepsilon,\ell})=+\infty$. Then
$\sup_{n\in\N} a_n  = +\infty$, where $a_n$ are partial quotients of the continued fraction of $\varepsilon$.
\end{prop}

\pfz
Since ${\rm ind}(u_{\varepsilon,\ell})=+\infty$, we can find for every $j\in\N$ a factor $w\in\L(u_{\varepsilon,\ell})$ such that $w^j\in\L(u_{\varepsilon,\ell})$. This means, using Theorem~\ref{thm:ABMP}, that $\sigma_{01}(w^j)=\big(\sigma_{01}(w)\big)^j\in\L(u_\varepsilon)$.
Any Sturmian word with slope $\varepsilon$ therefore contains arbitrarily long integer powers, and thus ${\rm ind}(u_\varepsilon)=+\infty$. The statement
of the proposition then follows from Theorem~\ref{thm:1}.
\pfk

Theorem~\ref{thm:3} implies that in case of a Sturmian word with finite index, the integer power occurring in $u_\varepsilon$ is at most
$2+\sup_{k\geq 1} a_k$. A direct consequence of the proof of Proposition~\ref{p:1} is the following corollary.

\begin{coro}
For integer exponents $j$ of powers of factors occurring in a 3iet word $u_{\varepsilon,\ell}$ we have
$$
j\leq 2+ \sup_{n\in\N} a_n\,.
$$
\end{coro}

For the demonstration of the implication opposite to Proposition~\ref{p:1} we use a more detailed description of the structure of Sturmian words,
as can be found in~\cite{lothaire}. For $\varepsilon=[0,a_1,a_2,a_3,\dots]$ we define a sequence of words $s_n$ recursively by
$$
s_{-1}:=1,\quad s_{0}:=0,\quad s_1:=s_0^{a_1-1}s_{-1},\quad\hbox{ and for $n\geq 2$}\quad s_{n+1}:=s_n^{a_{n+1}}s_{n-1}\,.
$$
It is known that the word defined by
$$
c_\varepsilon:=\lim_{n\to\infty} s_n
$$
is a Sturmian word with slope $\varepsilon$. The infinite word $c_\varepsilon$ plays an outstanding role among
Sturmian words of the same slope. Usually, it is called the characteristic word of $\varepsilon$. The languages $\L(c_\varepsilon)$
and $\L(u_\varepsilon)$ coincide for every Sturmian word $u_\varepsilon$ with slope $\varepsilon$. We will make use of the block structure
of the word $c_\varepsilon$, see~\cite{DaLenzNecelo}.

For arbitrary fixed $n\in\N$ we denote
\begin{equation}\label{eq:11}
k:=a_{n+1}\,,\quad E:=s_n\,,\quad F:=s_{n-1}\,.
\end{equation}
Then the infinite word $c_\varepsilon$ is a concatenation of two blocks, namely $E^{k+1}F$ and $E^kF$.
From the definition of the sequence $(s_n)$ it follows that for $n\geq 2$ we have $|E|>|F|$ and that  the first letters of the words $E$ and $F$ coincide.
In order to explain what this means for the structure of a Sturmian word $u_\varepsilon$, we recall the notion of {\em shift} $\Gamma:\A^\N\mapsto\A^\N$
on an infinite word $u=u_0u_1u_2\cdots$,
$$
\Gamma(u_0u_1u_2\cdots) = u_1u_2u_3\cdots\,.
$$
On the set of finite words we can define the {\em cyclic shift} $\gamma:\A^*\mapsto\A^*$ by
$$
\gamma(w_0w_1\cdots w_{n-1})=w_1\cdots w_{n-1}w_0\,.
$$
Since $u_\varepsilon$ has the same language as the characteristic word $c_\varepsilon$, the infinite word $u_\varepsilon$ is (up to a prefix)
also a concatenation of blocks $E^{k+1}F$, $E^kF$. Formally, for every $u_\varepsilon$ there exists an index $i_0$ such that the
infinite word $\Gamma^{i_0}(u_\varepsilon)$ is a concatenation of blocks $E^{k+1}F$ and $E^kF$.

\begin{prop}\label{p:2}
Let $u_{\varepsilon,\ell}$ be a 3iet word with parameters $\varepsilon,\ell$, and let the continued fraction $\varepsilon=[0,a_1,a_2,a_3,\dots]$
have bounded partial quotients. Then
$$
{\rm ind}(u_{\varepsilon,\ell})\geq  \sup_{n\in\N} \Big\lfloor \frac{a_n}{2}\Big\rfloor\,.
$$
\end{prop}

\pfz
Denote by $u^{(1)}$, $u^{(2)}$ the following Sturmian words with slope $\varepsilon$
$$
u^{(1)} := \sigma_{01} (u_{\varepsilon,\ell})\,,\qquad u^{(2)} := \sigma_{10} (u_{\varepsilon,\ell})\,.
$$
We know that $u^{(1)}$ is amicable with $u^{(2)}$. Thus for arbitrary iteration $i\in\N$ of the shift we have
$$
\Gamma^{i}(u^{(1)}) \propto \Gamma^{i}(u^{(2)}) \quad\hbox{ or }\quad
\Gamma^{i+1}(u^{(1)}) \propto \Gamma^{i+1}(u^{(2)})\,.
$$
Our aim is to find a factor $w$ in $u_{\varepsilon,\ell}$ which has in $u_{\varepsilon,\ell}$ sufficiently long power $w^m$.
We will construct the factor $w$ as the ternarization of pairs of amicable words $E$ and cyclic shifts of $E$.
We will use the notation of~\eqref{eq:11}. We distinguish two cases.

\begin{enumerate}

\item Suppose that there exists an iteration $i$ of the shift and $m\in\N$ such that
$$
\Gamma^i(u^{(1)}) = E_2E^mFr^{(1)}\,,\qquad \Gamma^i(u^{(2)}) = E^{k+1}Fr^{(2)}\,,
$$
where $r^{(1)}, r^{(2)}$ are suffixes of infinite words $u^{(1)}$, $u^{(2)}$ and $E_2$ is a suffix (possibly empty) of the factor $E$.
The situation is depicted on Figure~\ref{f}.

\begin{figure}[ht]
\setlength{\unitlength}{5.5pt}
\begin{center}
\begin{picture}(55,15)
\put(0,10){\makebox(0,0){$u^{(1)}=\quad\cdots$}}
\put(0,5){\makebox(0,0){$u^{(2)}=\quad\cdots$}}
\put(42,10){\makebox(0,0){$\cdots$}}
\put(59,5){\makebox(0,0){$\cdots$}}
\put(8,10){\line(1,0){30}}
\put(9,9.5){\line(0,1){1}}
\put(17,9.5){\line(0,1){1}}
\put(25,9.5){\line(0,1){1}}
\put(33,9.5){\line(0,1){1}}
\put(37,9.5){\line(0,1){1}}
\put(6,5){\line(1,0){50}}
\put(7,4.5){\line(0,1){1}}
\put(11,4.5){\line(0,1){1}}
\put(19,4.5){\line(0,1){1}}
\put(27,4.5){\line(0,1){1}}
\put(35,4.5){\line(0,1){1}}
\put(43,4.5){\line(0,1){1}}
\put(51,4.5){\line(0,1){1}}
\put(55,4.5){\line(0,1){1}}
\put(13,11.5){\makebox(0,0){$E$}}
\put(21,11.5){\makebox(0,0){$E$}}
\put(29,11.5){\makebox(0,0){$E$}}
\put(35,11.5){\makebox(0,0){$F$}}
\put(9,3.5){\makebox(0,0){$F$}}
\put(15,3.5){\makebox(0,0){$E$}}
\put(23,3.5){\makebox(0,0){$E$}}
\put(31,3.5){\makebox(0,0){$E$}}
\put(39,3.5){\makebox(0,0){$E$}}
\put(47,3.5){\makebox(0,0){$E$}}
\put(53,3.5){\makebox(0,0){$F$}}
\put(11,1.5){\vector(0,1){2.5}}
\put(11,0.5){\makebox(0,0){$i$}}
\put(11,13.0){\vector(0,-1){2.5}}
\put(11,14.0){\makebox(0,0){$i$}}
\end{picture}
\end{center}
\caption{}
\label{f}
\end{figure}

Then there exists a prefix $E_1$ of $E$ such that $E=E_1E_2$ and
$$
\Gamma^i(u^{(1)})=(E_2E_1)^mE_2Fr^{(1)}\,.
$$
Since $E$ and $F$ start with the same letter, we can write
$$
\Gamma^{i+1}(u^{(1)}) = \big(\gamma(E_2E_1)\big)^m\hat r^{(1)}\,,\qquad \Gamma^{i+1}(u^{(2)}) = \big(\gamma(E)\big)^{k+1}\hat r^{(2)}\,,
$$
where $\hat r^{(1)}, \hat r^{(2)}$ are some suffixes of infinite words $u^{(1)}$, $u^{(2)}$, respectively.
Since  $u^{(1)}\propto u^{(2)}$, we have either $\Gamma^{i}(u^{(1)})\propto \Gamma^{(i)}(u^{(2)})$ or
$\Gamma^{i+1}(u^{(1)})\propto \Gamma^{(i+1)}(u^{(2)})$. In the former case, $E_2E_1\propto E$, in the latter, $\gamma(E_2E_1)\propto \gamma(E)$. Due to the repetitions in $u^{(1)}$ and $u^{(2)}$ we can then infer (using Remark  \ref{uniqueness}) that
the infinite word $u_{\varepsilon,\ell}$ thus contains factors $\big({\rm ter}(E_2E_1, E)\big)^m$ or $\big({\rm ter}(\gamma(E_2E_1),\gamma(E))\big)^m$.

\item
Let us suppose that $i$ is an iteration of the shift such that
$$
\Gamma^i(u^{(1)}) = E^{k+1}Fr^{(1)}\,,\qquad \Gamma^i(u^{(2)}) =  (E_2E_1)^mFr^{(2)}\,,
$$
where $r^{(1)}, r^{(2)}$ are some suffixes of infinite words $u^{(1)}$, $u^{(2)}$, as depicted in Figure~\ref{ff}.
Analogously to the case 1, one finds the power $w^m$ in the infinite word $u_{\varepsilon,\ell}$.

\begin{figure}[ht]
\setlength{\unitlength}{5.5pt}
\begin{center}
\begin{picture}(55,15)
\put(0,10){\makebox(0,0){$u^{(1)}=\quad\cdots$}}
\put(0,5){\makebox(0,0){$u^{(2)}=\quad\cdots$}}
\put(42,5){\makebox(0,0){$\cdots$}}
\put(59,10){\makebox(0,0){$\cdots$}}
\put(8,5){\line(1,0){30}}
\put(9,4.5){\line(0,1){1}}
\put(17,4.5){\line(0,1){1}}
\put(25,4.5){\line(0,1){1}}
\put(33,4.5){\line(0,1){1}}
\put(37,4.5){\line(0,1){1}}
\put(6,10){\line(1,0){50}}
\put(7,9.5){\line(0,1){1}}
\put(11,9.5){\line(0,1){1}}
\put(19,9.5){\line(0,1){1}}
\put(27,9.5){\line(0,1){1}}
\put(35,9.5){\line(0,1){1}}
\put(43,9.5){\line(0,1){1}}
\put(51,9.5){\line(0,1){1}}
\put(55,9.5){\line(0,1){1}}
\put(13,3.5){\makebox(0,0){$E$}}
\put(21,3.5){\makebox(0,0){$E$}}
\put(29,3.5){\makebox(0,0){$E$}}
\put(35,3.5){\makebox(0,0){$F$}}
\put(9,11.5){\makebox(0,0){$F$}}
\put(15,11.5){\makebox(0,0){$E$}}
\put(23,11.5){\makebox(0,0){$E$}}
\put(31,11.5){\makebox(0,0){$E$}}
\put(39,11.5){\makebox(0,0){$E$}}
\put(47,11.5){\makebox(0,0){$E$}}
\put(53,11.5){\makebox(0,0){$F$}}
\put(11,2.0){\vector(0,1){2.5}}
\put(11,1.0){\makebox(0,0){$i$}}
\put(11,13.5){\vector(0,-1){2.5}}
\put(11,14.5){\makebox(0,0){$i$}}
\end{picture}
\end{center}
\caption{}
\label{ff}
\end{figure}

\end{enumerate}

Let us now find $m\geq \lfloor \frac{k}{2}\rfloor$ such that situation 1 or 2 occurs. Again, we discuss two cases.

\begin{itemize}
\item
Suppose we find the occurrences of the factor $F$ in the words $u^{(1)}$, $u^{(2)}$ such that indices occupied by $F$ in $u^{(1)}$
and indices occupied by $F$ in $u^{(2)}$ have a non-empty overlapping, see Figure~\ref{ffff}. In this case clearly $m\geq k-1$.

\begin{figure}[ht]
\setlength{\unitlength}{5.5pt}
\begin{center}
\begin{picture}(55,15)
\put(0,10){\makebox(0,0){$u^{(1)}=\quad\cdots$}}
\put(0,5){\makebox(0,0){$u^{(2)}=\qquad\quad\cdots$}}
\put(40.5,5){\makebox(0,0){$\cdots$}}
\put(38,10){\makebox(0,0){$\cdots$}}
\put(8,5){\line(1,0){30}}
\put(9,4.5){\line(0,1){1}}
\put(17,4.5){\line(0,1){1}}
\put(21,4.5){\line(0,1){1}}
\put(29,4.5){\line(0,1){1}}
\put(37,4.5){\line(0,1){1}}
\put(6,10){\line(1,0){30}}
\put(7,9.5){\line(0,1){1}}
\put(15,9.5){\line(0,1){1}}
\put(19,9.5){\line(0,1){1}}
\put(27,9.5){\line(0,1){1}}
\put(35,9.5){\line(0,1){1}}
%\put(43,9.5){\line(0,1){1}}
%\put(51,9.5){\line(0,1){1}}
%\put(55,9.5){\line(0,1){1}}
%
\put(13,3.5){\makebox(0,0){$E$}}
\put(19,3.5){\makebox(0,0){$F$}}
\put(25,3.5){\makebox(0,0){$E$}}
\put(33,3.5){\makebox(0,0){$E$}}
\put(11,11.5){\makebox(0,0){$E$}}
\put(17,11.5){\makebox(0,0){$F$}}
\put(23,11.5){\makebox(0,0){$E$}}
\put(31,11.5){\makebox(0,0){$E$}}
%\put(39,11.5){\makebox(0,0){$E$}}
%\put(47,11.5){\makebox(0,0){$E$}}
%\put(53,11.5){\makebox(0,0){$F$}}
%
\put(21,1.5){\vector(0,1){2.5}}
\put(21,0.5){\makebox(0,0){$i$}}
\put(21,13.0){\vector(0,-1){2.5}}
\put(21,14.0){\makebox(0,0){$i$}}
\end{picture}
\end{center}
\caption{}
\label{ffff}
\end{figure}

\item In the opposite case, consider iteration $t$ of the shift such that
$$
\Gamma^t(u^{(1)}) = E^{k+1}Fr^{(1)}
$$
for some suffix $r^{(1)}$ of $u^{(1)}$ and find the minimal index $j\geq 0$ such that $\Gamma^{t+j}(u^{(2)})=FE^kr^{(2)}$.
The block structure of $u^{(1)}$ ensures that $j<|E|(k+1)-|F|$, see Figure~\ref{fffff}.

\begin{figure}[ht]
\setlength{\unitlength}{5.5pt}
\begin{center}
\begin{picture}(55,16)
\put(0,14){\makebox(0,0){$u^{(1)}=\quad\cdots$}}
\put(0,9){\makebox(0,0){$u^{(2)}=\quad\cdots$}}
\put(65,9){\makebox(0,0){$\cdots$}}
\put(59,14){\makebox(0,0){$\cdots$}}
\put(8,9){\line(1,0){54}}
\put(9,8.5){\line(0,1){1}}
\put(17,8.5){\line(0,1){1}}
\put(25,8.5){\line(0,1){1}}
\put(33,8.5){\line(0,1){1}}
\put(37,8.5){\line(0,1){1}}
\put(45,8.5){\line(0,1){1}}
\put(53,8.5){\line(0,1){1}}
\put(61,8.5){\line(0,1){1}}
\put(6,14){\line(1,0){50}}
\put(7,13.5){\line(0,1){1}}
\put(11,13.5){\line(0,1){1}}
\put(19,13.5){\line(0,1){1}}
\put(27,13.5){\line(0,1){1}}
\put(35,13.5){\line(0,1){1}}
\put(43,13.5){\line(0,1){1}}
\put(51,13.5){\line(0,1){1}}
\put(55,13.5){\line(0,1){1}}
\put(13,7.5){\makebox(0,0){$E$}}
\put(21,7.5){\makebox(0,0){$E$}}
\put(29,7.5){\makebox(0,0){$E$}}
\put(35,7.5){\makebox(0,0){$F$}}
\put(41,7.5){\makebox(0,0){$E$}}
\put(49,7.5){\makebox(0,0){$E$}}
\put(57,7.5){\makebox(0,0){$E$}}
\put(9,15.5){\makebox(0,0){$F$}}
\put(15,15.5){\makebox(0,0){$E$}}
\put(23,15.5){\makebox(0,0){$E$}}
\put(31,15.5){\makebox(0,0){$E$}}
\put(39,15.5){\makebox(0,0){$E$}}
\put(47,15.5){\makebox(0,0){$E$}}
\put(53,15.5){\makebox(0,0){$F$}}
\put(11,4.5){\vector(0,1){4}}
\put(11,3.5){\makebox(0,0){$t$}}
%\put(11,13.5){\vector(0,-1){2.5}}
%\put(11,14.5){\makebox(0,0){$t$}}
%
\put(33,4.5){\vector(0,1){3.5}}
\put(33,3.5){\makebox(0,0){$t\!+\!j$}}
%\put(33,13.5){\vector(0,-1){2.5}}
%\put(33,14.5){\makebox(0,0){$t\!+\!j$}}
%
\put(37,1.5){\vector(0,1){6.5}}
\put(37,0.5){\makebox(0,0){$t\!+\!j\!+\!|F|$}}
%\put(37,14.5){\vector(0,-1){3.5}}
%\put(37,15.5){\makebox(0,0){$t\!+\!j\!+\!|F|$}}
%
\put(51,1.5){\vector(0,1){7.0}}
\put(51,0.5){\makebox(0,0){$t\!+\!(k\!+\!1)|E|$}}
\put(11.2,5.5){\vector(1,0){21.6}}
\put(32.8,5.5){\vector(-1,0){21.6}}
\put(22,4.5){\makebox(0,0){$x$}}
\put(37.2,4.0){\vector(1,0){13.6}}
\put(50.8,4.0){\vector(-1,0){13.6}}
\put(44,3.0){\makebox(0,0){$y$}}
\end{picture}
\end{center}
\caption{}
\label{fffff}
\end{figure}

As $x:=j$ and $y:=(k+1)|E|-j-|F|$ satisfy $x+y=(k+1)|E|-|F|$, either $x$ or $y$ is greater than $\frac12\big((k+1)|E|-|F|\big)$.
If the greater one is $x$, we put $i:=t$, otherwise $i:=t+j+|F|$. In both cases
$$
m\geq \Big\lfloor \frac{(k+1)|E|-|F|}{2|E|} \Big\rfloor \geq \Big\lfloor \frac{k}{2} \Big\rfloor\,.
$$
\end{itemize}
\pfk

%%%%%%%%%%%%%%%%%%%%%%%%%%%%%%%%%%%%%%%%%%%%%%%%%%%%%%%%%%%%%%%%%%%%%
\section{Comments}\label{sec:relationtoSturmanndRotation}

\begin{enumerate}
\item[(i)]
Theorems~\ref{thm:1} and~\ref{thm:XY} say that for arbitrary $\varepsilon$, any 3iet word $u_{\varepsilon,\ell}$ and Sturmian word $u_\varepsilon$
satisfy
$$
{\rm ind}(u_{\varepsilon,\ell})<+\infty \quad\iff\quad {\rm ind}(u_{\varepsilon})<+\infty\,.
$$
Boundedness of the index of a 3iet word $u_{\varepsilon,\ell}$ thus does not depend on the parameter $\ell$.

\item[(ii)]
If ${\rm ind}(u_{\varepsilon,\ell})$ is finite, we have an estimate
\begin{equation}\label{eq:bound}
\Big\lfloor\frac{K}{2}\Big\rfloor \leq {\rm ind}(u_{\varepsilon,\ell}) \leq K+3\,,
\end{equation}
where $K:=\max_{n\in\N} a_n$. One can ask, whether the index of a 3iet word can be close to both the upper and the lower bounds, depending  on the parameter $\ell$. Recall that $\ell$ satisfies $\max\{\varepsilon,1-\varepsilon\}<\ell<1$.

Consider a sequence $\ell_n\in(0,1)$ with $\lim_{n\to\infty}\ell_n=1$, and a sequence of transformations $T_{\varepsilon,\ell_n}$. Let us denote by $u^{(n)}$ the 3iet word coding the orbit of $x_0=0$ under the transformation $T_{\varepsilon,\ell_n}$. The length of the interval $I_B:=[\ell_n-1+\varepsilon,\varepsilon)$ is approaching 0, as $n$ tends to infinity, therefore the density of the letter $B$ in words $u^{(n)}$
also decreases to 0. Obviously,
$$
\lim_{n\to\infty} u^{(n)} = u_\varepsilon\,,
$$
where $u_\varepsilon$ is the Sturmian word in the alphabet $\{A,C\}$ coding the orbit of 0 under $T_\varepsilon$.
Since $u_\varepsilon$ contains integer powers $K+1$, the index of infinite words $u^{(n)}$ for sufficiently large $n$ satisfies
$$
{\rm ind}(u^{(n)})\geq K+1.
$$
The index of 3iet words will take values greater than $K+1$ also if we consider parameters $\ell$ approaching to $\max\{\varepsilon,1-\varepsilon\}$.
The limit results in a Sturmian word in the alphabet $\{A,B\}$ or $\{B,C\}$.

How well the lower bound in~\eqref{eq:bound} can be approached, remains an open question.

\item[(iii)]
As mentioned in Section~\ref{sec:preliminaries}, every word $v$ coding a rotation is a morphic image of a 3iet word $u$. We show that ${\rm ind}(v)$ and ${\rm ind}(u)$ can substantially differ. Consider the sequence of words $u^{(n)}$ defined in (ii) with slope $\varepsilon$ with bounded partial quotients.
Take the morphism $\Phi_0: A\mapsto 0, B\mapsto 01, C\mapsto 0$. The index of $u^{(n)}$ is then bounded by $K+3$, whereas the words $\Phi_0(u^{(n)})$ coding rotations have as  limit the periodic word $0000\cdots$, and thus ${\rm ind}(\Phi_0(u^{(n)}))\to+\infty$.

\item[(iv)] For Sturmian words with parameter $\varepsilon$ it is known  \cite{Mignosi,Durand} that the three properties, namely finite index, boundedness of the coefficients of the  continued fraction expansion of $\varepsilon$, and linear recurrence are all equivalent. This is not the case anymore for 3iets as discussed in  \cite{ferenczi-holton-zamboni-join-etds-25} based on some (unpublished) work of Boshernitzan.

\item[(v)]  Powers in  words play a role in the study of the associated Schroedinger operators. In fact, due to the so-called Gordon argument, occurrence of sufficiently many powers of order three (and higher) can be used to exclude eigenvalues (see \cite{damanik} for a survey on this topic and further references). Whenever the continued fraction expansion of $\varepsilon$ has infinitely many coefficients with value at least six, we can then conclude (almost sure) absence of eigenvalues.  This gives a somewhat more explicit version of a result of \cite{CGdO}. Note, however, that the considerations of \cite{CGdO} are not confined to 3iets.

\end{enumerate}

%%%%%%%%%%%%%%%%%%%%%%%%%%%%%%%%%%%%%%%%%%%%%%%%%%%%%%%%%%%%%%%%%%%%%
\section*{Acknowledgements}

We acknowledge financial support by the Czech Science Foundation
grant 201/09/0584 and by the grants MSM6840770039 and LC06002 of
the Ministry of Education, Youth, and Sports of the Czech
Republic.

%%%%%%%%%%%%%%%%%%%%%%%%%%%%%%%%%%%%%%%%%%%%%%%%%%%%%%%%%%%%%%%%%%%%%

%%%%%%%%%%%%%%%%%%%%%%%%%%%%%%%%%%%%%%%%%%%%%%%%%%%%%%%%%%%%%%%%%%%%%
\end{document}